\documentclass[a4paper,11pt]{amsart}
\usepackage{MnSymbol}
\newcommand{\x}{\operatorname\ell}

\newcommand{\Li}{\operatorname{\mathcal L}}

\begin{document}
\title{Pentagon Identity Revisited I}
\author[A. Yu. Volkov]{Alexander Yu. Volkov}
\address{Steklov Mathematical Institute, St. Petersburg}
\thanks{This work has been supported by the Swiss NSF (grant 200020--126817)}
\begin{abstract}
This is to present the previously overlooked q-algebraic relation
\begin{equation*}
 \x(x)\x(y)=\x\left((1-x)^{-1}y\right)
 \x\left(-x(1-x-y)^{-1}y\right)\x\left(x(1-y)^{-1}\right),
\end{equation*}
where \(x\) and \(y\) are q-commuting variables, \(yx=qxy\), and \(\x\) is Heine's q-exponential function, \(\x(z)=1+z/(1-q)+z^2/(1-q^2)(1-q)+\ldots\,\). This relation represents the most direct q-analog of the classical five-term dilogarithm relation
\begin{equation*}
 \Li(x)+\Li(y)
 =\Li\left(\frac y{1-x}\right)
 +\Li\left(\frac{-xy}{1-x-y}\right)
 +\Li\left(\frac x{1-y}\right),
\end{equation*}
where \(x\) and \(y\) are ordinary commuting variables, \(\Li(z)=L(1/(1-1/z))\), and \(L\) is the Rogers dilogarithm.
\end{abstract}
\maketitle
\section{Sch\"utzenberger's Formula}
Let \(x\) and \(y\) be ``q-commuting variables'', that is elements of a suitable noncommutative algebra satisfying Weyl's relation
\begin{equation*}
 yx=qxy,
\end{equation*}
with \(q\) being a constant (not equal 1).

Clearly, this type of noncommutativity leads to the same binomial formula as in the classical case,
\begin{equation*}
 (x+y)^n=\sum_{k=0}^nc_k^n\,x^{n-k}y^k,
\end{equation*}
albeit with different coefficients. To find them, note that
\begin{multline*}
 (x+y)^n=(x+y)^{n-1}(x+y)=\sum c_k^{n-1}x^{n-k-1}y^k(x+y)\\
 =\sum c_k^{n-1}(q^kx^{n-k}y^k+x^{n-k-1}y^{k+1})
 =\sum (q^kc_k^{n-1}+c_{k-1}^{n-1})x^{n-k}y^k,
\end{multline*}
and hence the q-Pascal triangle equation has the form
\(c_k^n=q^kc_k^{n-1}+c_{k-1}^{n-1}\).
The solution is then easily guessed to be
\begin{equation*}
 c_k^n=\frac{(1-q^n)(1-q^{n-1})\ldots(1-q^{n-k+1})}
 {(1-q^k)(1-q^{k-1})\ldots(1-q)}
 =\frac{(q)_n}{(q)_{n-k}(q)_k},
\end{equation*}
where the usual notation
\((q)_n=(1-q^n)(1-q^{n-1})\ldots(1-q)\) is employed.

Now, rewrite the binomial formula as
\begin{equation*}
 \frac{(x+y)^n}{(q)_n}=\sum_{k=0}^n
 \frac{x^{n-k}}{(q)_{n-k}}\frac{y^k}{(q)_k}
\end{equation*}
and then (formally) sum it over \(n\geq0\) to obtain Sch\"utzenberger's formula~\cite{s}
\begin{equation*}
 \x(x+y)=\x(x)\x(y),
\end{equation*}
where \(\x\) is the now ubiquitous ``q-exponential function''
\begin{equation*}
 \x(z)=\sum_{n=0}^\infty\frac{z^n}{(q)_n}.
\end{equation*}

In view of this discovery, it indeed seems apt to call this series the q-exponential function --- yet it has also become known as ``the quantum dilogarithm''. We will soon see why.

\section{Pentagon Identity}
Let us compute the reverse product \(\x(y)\x(x)\). We have
\begin{multline*}
 \x(y)\x(x)=\x(y)\x(x)
 \underbrace{\x(x)^{-1}\x(x+y)\x(y)^{-1}}_
 {\text{=1 by Sch\"utzenberger's formula}\mspace{-260mu}}\\
 =\x(y)\x(x+y)\x(y)^{-1}
 =\x\left(\x(y)x\x(y)^{-1}+y\right),
\end{multline*}
and thus it boils down to finding the conjugation of \(x\) by \(\x(y)\). To this end, note that \((q)_n=(1-q^n)(q)_{n-1}\), and therefore
\begin{equation*}
 \x(z)-\x(qz)=\sum_{n=0}^\infty\frac{(1-q^n)z^n}{(q)_n}
 =\sum_{n=1}^\infty\frac{z^n}{(q)_{n-1}}
 =z\x(z),
\end{equation*}
or, in a more transparent form,
\begin{equation*}
 \frac{\x(qz)}{\x(z)}=1-z.
\end{equation*}
Hence
\begin{equation*}
 \x(y)x\x(y)^{-1}=x\x(qy)\x(y)^{-1}=x(1-y),
\end{equation*}
and thus, finally,
\begin{equation*}
 \x(y)\x(x)=\x(x+y-xy).
\end{equation*}

This relation (due to Faddeev and myself~\cite{fv}) is still a bit of a mystery, but it has at least one useful consequence. Following Faddeev and Kashaev~\cite{fk}, note that \(x-xy\) and \(y\) satisfy the same relation as \(x\) and \(y\):
\(y(x-xy)=q(x-xy)y\). Hence, by Sch\"utzenberger's formula,
\(\x(x+y-xy)=\x(x-xy)\x(y)\). For the same reason, \(\x(x-xy)=\x(x)\x(-xy)\), and thus, finally,
\begin{equation*}
 \x(y)\x(x)=\x(x)\x(-xy)\x(y).
\end{equation*}

This relation is now widely known as the ``pentagon identity'', and its popularity is partly due to the following reason. According to Faddeev and Kashaev, if \(q\) is sent to 1 in a certain nontrivial way, this five-factor relation degenerates into the classical five-term dilogarithmic relation
\begin{equation*}
 \Li(x)+\Li(y)
 =\Li\left(\frac y{1-x}\right)
 +\Li\left(\frac{-xy}{1-x-y}\right)
 +\Li\left(\frac x{1-y}\right),
\end{equation*}
where \(\Li(z)=L(1/(1-1/z))\), and \(L\) is the Rogers dilogarithm.

This observation is remarkable in more ways than one. On the one hand, it shows that in the q-world the same function can indeed be both exponential and dilogarithmic at the same time. On the other hand, it alarms us that something is amiss here. One normally expects closer similarity between q-formulas and their classical counterparts, so there is clearly a missing link between the above two relations. We will now find out what it is. 
\section{True Five-factor Relation}
The computation here is quite similar to the one in the beginning of Section~2. We have
\begin{multline*}
 \x(x)\x(y)=\x(x)\x(y)
 \underbrace{\x(x)^{-1}\x(y)^{-1}\x(x)\x(-xy)\x(y)}_
 {\text{=1 by the pentagon identity}\mspace{-235mu}}\\
 =\x(x)\x(y)\x(x)^{-1}\x(y)^{-1}\x(x)\x(-xy)
 \underbrace{\x(x)^{-1}\x(y)\x(y)^{-1}\x(x)}_
 {\text{=1 by itself}\mspace{-105mu}}
 \x(y)\\
 =\x\left(\x(x)y\x(x)^{-1}\right)
 \x\left(-\x(y)^{-1}\x(x)xy\x(x)^{-1}\x(y)\right)
 \x\left(\x(y)^{-1}x\x(y)\right),
\end{multline*}
and this time it boils down to computing the conjugation of \(x\), \(y\) and \(xy\) by \(\x(y)^{-1}\), \(\x(x)\) and \(\x(y)^{-1}\x(x)\) respectively. The first two are computed in exactly the same way as the one in Section~2,
\begin{align*}
 \x(y)^{-1}x\x(y)&=x\x(qy)^{-1}\x(y)=x(1-y)^{-1}\\
 \x(x)y\x(x)^{-1}&=\x(x)\x(qx)^{-1}y=(1-x)^{-1}y,
\end{align*}
while the third one is a combination of these:
\begin{multline*}
 \x(y)^{-1}\x(x)xy\x(x)^{-1}\x(y)
 =\x(y)^{-1}x(1-x)^{-1}y\x(y)\\
 =x(1-y)^{-1}(1-x(1-y)^{-1})^{-1}y
 =x(1-x-y)^{-1}y.
\end{multline*}
Finally, putting it all together gives the announced five-factor relation:
\begin{equation*}
 \x(x)\x(y)=\x\left((1-x)^{-1}y\right)
 \x\left(-x(1-x-y)^{-1}y\right)\x\left(x(1-y)^{-1}\right).
\end{equation*}

Our immediate goal is thus achieved --- indeed, one can hardly imagine a truer q-analog of the five-term dilogarithm relation (at the bottom of the previous page). Still, it would perhaps be interesting to analyze the precise limiting procedure leading from one relation to the other; this work is currently in progress~\cite{kn}.

Also, a few words are in order about how this newfound five-factor relation fits into some larger context.

\section{\(\mathrm{qY}\)-system}
Recall that the five arguments in the five-term relation are related to one another by a so-called Y-system. Indeed, if we denote by
\(Y_t,t\in\mathbb Z/5\mathbb Z\) the following five expressions in (commuting) variables \(x\) and \(y\):
\begin{align*}
 Y_1&=y&Y_2&=\frac{1-y}x&Y_3&=\frac{1-x-y}{-xy}&Y_4&=\frac{1-x}y&Y_5&=x,
\end{align*}
then the five-term relation takes the form
\begin{equation*}
 \Li(Y_5)+\Li(Y_1)
 =\Li(1/Y_4)+\Li(1/Y_3)+\Li(1/Y_2),
\end{equation*}
and these \(Y\)'s are easily checked to satisfy the five equations
\begin{equation*}
 Y_{t-1}Y_{t+1}=1-Y_t,
\end{equation*}
which are easily recognized (up to the overall sign) as the Y-system of type \((A_1,A_2)\) (see \cite{k} and references therein).

Obviously, the q-case lends itself to similar treatment --- except here we have two relations to deal with. We write the pentagon identity of Section~3 as
\begin{equation*}
 \x\left(X_1^{}\right)\x\left(X_5^{}\right)
 =\x\left(qX_2^{-1}\right)
 \x\left(qX_3^{-1}\right)\x\left(qX_4^{-1}\right),
\end{equation*}
and the five-factor relation of Section~4 as
\begin{equation*}
 \x\left(Y_5^{}\right)\x\left(Y_1^{}\right)
 =\x\left(qY_4^{-1}\right)
 \x\left(qY_3^{-1}\right)\x\left(qY_2^{-1}\right),
\end{equation*}
where the (cyclic) quintuples \(X_t\) and \(Y_t\), \(t\in\mathbb Z/5\mathbb Z\) are as follows:
\begin{align*}
 X&=\{y,qx^{-1},-q^2x^{-1}y^{-1},qy^{-1},x\}\\
 Y&=\{y,x^{-1}(q-y),-qx^{-1}(q-x-y)y^{-1},(q-x)y^{-1},x\}.
\end{align*}
These are easily checked to satisfy, on the one hand, the commutation relations
\(X_{t+1}X_t=qX_tX_{t+1}\)
and
\(Y_{t+1}Y_t=qY_tY_{t+1}\),
and on the other, the equations relating successive \(X\)'s and \(Y\)'s to one another:
\begin{equation*}
 X_{t-1}X_{t+1}=
 \begin{cases}
  \;\;q&\text{for \(t\equiv1,5\!\!\!\!\pmod5\)}\\
  -X_t&\text{for \(t\equiv2,3,4\!\!\!\!\pmod5\).}
 \end{cases}
\end{equation*}
and, for all \(t\),
\begin{equation*}
 Y_{t-1}Y_{t+1}=q-Y_t.
\end{equation*}

Of course, the latter equations should be taken as the definition of the ``qY-system of type \((A_1,A_2)\)'' and used as a starting point for a systematic theory of quantum Y-systems --- but it is perhaps more important that this qY-system has emerged here in a curious mix with two simpler ``qX-systems'',
\(X_{t-1}X_{t+1}=q\)
and
\(X_{t-1}X_{t+1}=-X_t\).
The relation of these to the qY-system is quite obvious: they are obtained by dropping one or the other term in the latter's right hand side. It is also easily checked that, taken separately, they have periods four and six --- in the sense that they are satisfied by the quadruple
\(\{X_1,X_2,X_4,X_5\}\)
and sextuple
\(\{X_1,X_2,X_3,X_4,X_5,qX_3^{-1}\}\)
respectively. This, in turn, must have something to do with the division of the five factors into two on the left hand side and three on the right hand side, so some explanation is in order as to how all these various observations might fit together. This will be given in the forthcoming Part~II of these notes.

To conclude, I thank Anton Alexeev, Ludwig Faddeev, Rinat Kashaev, Tomoki Nakanishi and Andras Szenes for numerous stimulating discussions.


\begin{thebibliography}{9}
\bibitem{s}Marcel~Paul Sch{\"u}tzenberger, \emph{Une interpr\'etation de certaines solutions de l'\'equation fonctionnelle: \(F(x+y)=F(x)F(y)\)}, C. R. Acad. Sci. Paris \textbf{236} (1953), 352--353.
\bibitem{fv}L.~Faddeev and A.~Yu.~Volkov, \emph{Abelian current algebra and the Virasoro algebra on the lattice}, Phys. Lett. B \textbf{315} (1993), no.~3-4, 311--318.
\bibitem{fk}L.~D.~Faddeev and R.~M.~Kashaev, \emph{Quantum dilogarithm}, Modern Phys. Lett. \textbf{A9} (1994), no.~5, 427--434.
\bibitem{kn}R.~Kashaev and T.~Nakanishi, in preparation.
\bibitem{k}B.~Keller, \emph{The periodicity conjecture for pairs of Dynkin diagrams}, arXiv:1001.1531.
\end{thebibliography}
\end{document}